\documentclass[12pt,leqno]{amsart}
\usepackage{amsmath,amstext,amssymb,amsopn,amsthm}

\title[stable process in parabola]
{Symmetric stable processes in parabola--shaped regions}
\author[R{.} Ba\~nuelos and K{.} Bogdan]
{Rodrigo Ba\~nuelos  and Krzysztof Bogdan}
\address{Department of Mathematics,
Purdue University,
West Lafayette, IN 47907-1395\newline
Institute of Mathematics, Polish Academy of Sciences,
and Institute of Mathematics, Wroc{\l}aw University of Technology,
50-370 Wroc{\l}aw, Poland}
\email{banuelos@math.purdue.edu  bogdan@im.pwr.wroc.pl}
\theoremstyle{plain}
\newtheorem{thm}{Theorem}
\newtheorem{lem}[thm]{Lemma}

\newcommand{\R}{\mathbb{R}}
\newcommand{\Rd}{{\mathbb{R}^d}}
\newcommand{\Rdl}{{\mathbb{R}^{d-1}}}
\newcommand{\tx}{\tilde{x}}
\newcommand{\ty}{\tilde{y}}
\newcommand{\tz}{\tilde{z}}
\newcommand{\tb}{{\tau_\beta}}
\newcommand{\ts}{{\tau_s}}
\newcommand{\Pb}{{\mathcal{P}_\beta}}
\newcommand{\PP}{{\bf P}}
\newcommand{\EE}{{\bf E}}
\DeclareMathOperator{\dist}{dist}

\theoremstyle{definition}

\theoremstyle{remark}
\newtheorem*{rem*}{Remark}

\date{7/14/2004}
\begin{document}

\sloppy
\footnotetext{
\emph{The first author was supported in part by NSF grant
\# 9700585-DMS}\\
 \emph{The second author was supported in part by KBN (2P03A 041 22)
 and by RTN (HPRN-CT-2001-00273-HARP)}
}
\begin{abstract}
We identify the critical exponent of integrability of the first exit time 
of rotation invariant stable L\'evy process from parabola--shaped region. 
\end{abstract}
\maketitle

\section{Introduction}
\label{sec:introduction}

 For $d=2,3,\ldots$ and $0<\beta<1$, we 
 define the {\it parabola--shaped
  region} in $\R^d$
$$
\Pb =\{x=(x_1,\tilde{x})\;:\; x_1>0\,,\; \tilde{x}\in \Rdl\,,|\tilde{x}|<x_1^\beta\}.
$$
Let $0<\alpha<2$.
By $\{X_t\}$ we denote the isotropic $\alpha$-stable $\Rd$-valued
L\'evy process (\cite{S}).  
The process is time-homogeneous, 
has right-continuous trajectories with left limits, $\alpha$-stable rotation
invariant independent increments, and characteristic  function
\begin{equation}
  \label{wfc}
\EE_xe^{i\xi (X_t-x)} = e^{-t |\xi|^\alpha}, \quad x\in \Rd\,,\quad
  \xi \in \Rd,  \quad  t \ge 0.
\end{equation}
Here $\EE_x$ is the expectation with respect to the
distribution $\PP_x$ of the process starting from $x \in \Rd$.
For an open set $U\subset \Rd$, we define  $\tau_U=\inf \{t\geq 0;\;
X_t\notin U\}$, the first exit time of $U$ (\cite{S}).  
In the case of the parabola-shaped region $\Pb$, we simply write 
$\tb$ for $\tau_\Pb$.  

The main result of this note is the following result.
\begin{thm}\label{iet} 
Let $p\geq 0$. Then 
$E_x \tb^p<\infty$ for (some, hence for all) $x\in \Pb$ if and only if
$p<p_0$, where 
$$p_0=\frac{(d-1)(1-\beta)+\alpha}{\alpha\beta}.$$
\end{thm} 

Theorem~\ref{iet} may be regarded as an addition to the research direction initiated
in \cite{BDS}, where it was proved that for $\beta=1/2$, $d=2$ and $\{X_t\}$
replaced by the {\it Brownian motion process} $\{B_{t}\}$, 
 $\tau_\beta$ is subexponentially integrable.  
For $\{B_t\}$, the generalizations to all the considered domains
$\Pb$ (along with essential strengthening of the result of
\cite{BDS}) were subsequently obtained in \cite{Li}, \cite{LS}, 
\cite{van}, \cite{DS} and \cite{BC}.  
We should note that this direction of research was influenced to a
large degree by the result of Burkholder \cite{Bk} on the critical
order of integrability of the exit times of $\{B_t\}$ from cones.  
For more on the many generalizations of Burkholder's result we refer the reader
to \cite{BS} and references therein. 

Brownian motion is a limiting case of the isotropic
$\alpha$-stable process and $B_{2t}$ corresponds to $\alpha=2$ via the
analogue of (\ref{wfc}). 
Extension of some of the above-mentioned results pertaining to cones
to the case $0<\alpha<2$ (that is, for $\{X_t\}$) were given in \cite{De2},
\cite{K3}, \cite{M} and \cite{BB}, see also \cite{BJ}. 
It should be emphasized that while there are many similarities between
$\{X_t\}$ and $\{B_{2t}\}$, there also exist essential differences in their respective 
properties and their proofs.
For example the critical exponent of integrability
of the exit time of $\{X_t\}$ is less than $1$ for every cone, however
narrow it may be (\cite{BB}), while it is arbitrarily large for 
$\{B_t\}$ in sufficiently narrow cones (\cite{Bk}).
A similar remark applies to regions $\Pb$.  The critical exponent  of integrability
of $\tau_\beta$ for $\{X_t\}$ given in Theorem~\ref{iet}
is {\it qualitatively different} from that of $\{B_t\}$ (\cite{BC}). 

Finally, we have recently
learned that Pedro J.
M\'endez-Hern\'andez, in a paper in preparation, has obtained results similar to those
presented.  His results apply to more general  regions defined by other ``slowly" increasing functions.

The proof of Theorem~\ref{iet} is based on estimates for harmonic
measure of $\{X_t\}$, following  the general idea used in  
\cite{BC} for $\{B_t\}$. 
The necessity of the condition $p<p_0$ for finiteness of
$E_x\tb^p$ is proved in Lemma~\ref{lA} and the sufficiency is proved in Lemma~\ref{lB}.
The main technical result of the paper is Lemma~\ref{lomh} where we
prove sharp estimates for harmonic measure of cut-offs of the
parabola-shaped region. 
At the end of the paper we discuss some remaining open problems.

\section{Proofs}
\label{sec:preliminaries}
We begin by reviewing the notation.
By $|\cdot|$ we denote the Euclidean norm in $\Rd$.
For $x\in \Rd$, $r>0$,
and a set $A \subset\Rd$ we let $B(x,r)=\{y \in \Rd \colon |x-y|<r\}$ and
$\dist(A,x)=\inf\{|x-y|\colon y \in A \}$. 
$A^{c}$ is the complement of $A$.
We always assume Borel measurability of the considered sets and functions.
In equalities and inequalities, unless stated otherwise,
$c$ denotes some unspecified positive real number whose value depends 
{\it only} on $d$, $\alpha$ and $\beta$.
 
It is well known that $(X_t,\PP_x)$ is a strong Markov processes with respect to the 
standard filtration (\cite{S}).

For an open set $U\subset \Rd$, with exit time $\tau_U$,  
the $\PP_x$ distribution of $X_{\tau_U}$:
$$
\omega^x_U(A):=P_x\{X_{\tau_U}\in A, \tau_U<+\infty\}\,,\quad A\subset \Rd\,,
$$ 
is a subprobability measure concentrated on $U^c$ (probability measure if $U$
is bounded) called the {\it harmonic measure}.
When $r>0$, $|x|<r$ and $B=B(0,r) \subset \Rd$, the corresponding harmonic measure
has the density function $P_r(x,y)=d \omega^x_B/dy$ (the \emph{Poisson
kernel} for the ball). We have
\begin{equation}
\label{Poisson}
   P_r(x,y) =
         C(d,\alpha)\left[ \frac{r^2-|x|^2}{|y|^2-r^2} \right] ^
{\alpha/2}|y-x|^{-d}
   \quad \text{if} \quad |y| >r,
\end{equation}
where $C(d,\alpha)= \Gamma(d/2)\pi^{-d/2-1}\sin(\pi\alpha/2)$, and  $P_r(x,y)=0$
otherwise (\cite{BGR}). 

The {\it scaling} property of $X_t$ plays a role in this paper.  
Namely, for every $r>0$ and $x\in \Rd$ the
$\PP_x$ distribution of $\{X_t\}$ is the same as the $\PP_{rx}$
distribution of $\{r^{-1}X_{r^{\alpha}t}\}$ (see (\ref{wfc})). 
In particular, the $\PP_x$ distribution of $\tau_{U}$ is the same as 
the $\PP_{rx}$ distribution of $r^{-\alpha} \tau_{rU}$. 
In short, $\tau_{rU}=r^\alpha \tau_U$ in distribution.

Let $\Pb'=\Pb\cap \{x_1>1\,,|\tilde{x}|<x_1^\beta/2\}$.
We claim that if $x\in \Pb'$ then $B(x,|x|^\beta/5)\subset \Pb$.
Indeed, let $x=(x_1,\tx)\in \Pb'$ and $y=(y_1,\ty)\in
B(x,|x|^\beta/5)$.
We have $|\tx|<x_1^\beta/2<x_1/2$ hence
$$
1<x_1\leq |x|=\sqrt{x_1^2+|\tx|^2}<\sqrt{5/4}x_1<5x_1/4\,.
$$
Then 
$$
y_1\geq x_1-|x|^\beta/5>x_1-|x|/5>x_1-(5/4)x_1/5=3x_1/4\,,
$$
and so
\begin{displaymath}
  |\ty| \leq |\tx|+|x|^\beta/5<x_1^\beta/2+(5x_1/4)^\beta/5
<3x_1^\beta/4<3(4y_1/3)^\beta/4<y_1^\beta\,,
\end{displaymath}
which yields $y\in \Pb$, as claimed.

\begin{lem}\label{lA}
  If $p\geq p_0$ then $E_x \tb^p=\infty$ for every $x\in \Pb$.
\end{lem}
\proof
Define $\tau=\inf\{t\geq 0:  X_t\notin B(X_0, |X_0|^\beta/5)\}$.
For $y\in \Rd$ and (nonnegative) $p$ by scaling we have
$
\EE_y\tau^p=\EE_0\tau_{B(0,1)}^p
(|y|^\beta/5)^{p\alpha}=c|y|^{\alpha\beta p}\,.
$
Let $x\in \Pb$, $r=\dist(x,\Pb^c)$, $B=B(x,r)$, and let $R$ be so
large that $B\subset \Pb\cap \{y_1\leq R\}$.
By strong Markov property we have
\begin{eqnarray*}
  \EE_x \tb^p
&\geq& 
\EE_x\{ X_{\tau_B}\in \Pb'\,;\;\EE_{X_{\tau_B}}\tb^p\}
\geq \EE_x\{X_{\tau_B}\in \Pb'\,;\; \EE_{X_{\tau_B}}\tau^p\}\\
&=&
c \EE_x\{X_{\tau_B}\in \Pb'\,;\; |X_{\tau_B}|^{\alpha\beta p}\}
\geq c \int_{\Pb'\cap \{y_1>R\}}P_r(0,y-x)|y|^{\alpha\beta p}\,dy\\
&\geq& cr^\alpha C(d,\alpha) \int_R^\infty dy_1\int_{|\ty|<y_1^\beta/2}
|y-x|^{-d-\alpha}|y|^{\alpha\beta p}\,d\ty\\
&\geq& const.\int_R^\infty y_1^{\beta(d-1)+\alpha\beta p-d-\alpha}\,dy_1\,.
\end{eqnarray*}
If $p\geq p_0$, then
$\beta(d-1)+\alpha\beta p-d-\alpha\geq
\beta(d-1)+ (d-1)(1-\beta)+\alpha-d-\alpha=-1$ and 
 the last
integral is positive infinity, proving the lemma.  
\qed

For $0\leq u<v\leq\infty$, we let $\Pb^{u,v}=\Pb\cap \{(x_1,\tx)\,:\;
u\leq x_1 <v\}$.
Let $0<s<\infty$ and define $\tau_s=\tau(s)=\tau_{\Pb^{0,s}}$.
Consider the cylinder 
\begin{equation}
  \label{edc}
\mathcal{C}=\{x=(x_1,\tx)\in \Rd\,:\; -\infty<x_1<s\,,\;
|\tx|<s^\beta\}\,.
\end{equation}
Clearly for $A\subset \Pb$, 
\begin{equation}
  \label{enh}
\PP_x\{X_\ts\in A\}\leq \PP_x\{X_{\tau_{\mathcal{C}}}\in
A\}\,,\quad x\in \Rd\,.
\end{equation}
By Lemma~4.3 of \cite{BK} for $A\subset\Rd\cap\{z_1\geq s\}$ we have
\begin{equation}
  \label{ewbk}
\PP_x\{X_{\tau_{\mathcal{C}}}\in
A\}\leq
c\int_{z=(z_1,\tz)\in A} \frac{s^{\alpha\beta}}{|z-x|^{d+\alpha}}
\left[
\frac{s^{\alpha\beta/2}}{(z_1-s)^{\alpha/2}}\vee 1\right]\,dz
\,,\quad x\in \mathcal{C}.
\end{equation}
The following lemma will simplify the use of the estimate (\ref{ewbk}).
\begin{lem}\label{ler}
  Let $s\geq 1$ and $x=(x_1,\tx)\in \Pb$ with $x_1\leq s/2$.
Let $s\leq u<v\leq \infty$ and assume that either $u\geq s+s^\beta$ or
$u=s$ and $v\geq s+s^\beta$.
Then
\begin{equation}
  \label{eer}
  \PP_x\{X_{\tau_s}\in \Pb^{u,v}\}
\leq
c s^{\alpha\beta}\int_u^v t^{-\alpha\beta p_0 -1} dt\,.
\end{equation}
\end{lem}
\proof
Denote
$$
f(y)=\frac{d \omega^x_{\Pb^{0,s}}}{dy}(y)\,.
$$
Let $y=(y_1,\ty)\in \Pb^{s,\infty}$. 
From (\ref{enh}) and (\ref{ewbk}) we can conclude that
$$
f(y)\leq c\frac{s^{\alpha\beta}}{|y-x|^{d+\alpha}}
\left[
\frac{s^{\alpha\beta/2}}{(y_1-s)^{\alpha/2}}\vee 1\right]\,.
$$
Since $s\geq 1$, we have $|\ty|\leq
y_1$ and $|y|\leq \sqrt{2}y_1$.
As $|y-x|\geq y_1-x_1\geq y_1/2$,
$$
f(y)\leq cs^{\alpha\beta}y_1^{-d-\alpha}
\left[\frac{s^{\alpha\beta/2}}{(y_1-s)^{\alpha/2}}\vee 1\right]\,.
$$
If $u\geq s+s^\beta$ then
\begin{eqnarray*}
\PP_x\{X_{\tau_s}\in \Pb^{u,v}\}
&\leq&
c s^{\alpha\beta}\int_u^v
y_1^{-d-\alpha}\int_{|\ty|<y_1^\beta}d\ty\,dy_1\\
&=&c s^{\alpha\beta}\int_u^v t^{-d-\alpha+\beta(d-1)}dt
=c s^{\alpha\beta}\int_u^v t^{-\alpha\beta p_0 -1} dt\,.
\end{eqnarray*}
If  $u=s$ and $v\geq s+s^\beta$, we consider
$$
I(s)=s^{\alpha\beta}\int_s^{s+s^\beta} t^{-\alpha\beta p_0 -1} dt\,,
$$
and
$$
II(s)=s^{\alpha\beta}\int_s^{s+s^\beta} t^{-\alpha\beta p_0 -1}
\frac{s^{\alpha\beta/2}}{(t-s)^{\alpha/2}} dt\,.
$$
It is enough to verify that $II(s)\leq cI(s)$.
Note that $s+s^\beta\leq 2s$, hence
$$
I(s)\geq s^{\alpha\beta} (2s)^{-\alpha\beta p_0 -1}s^\beta
= c s^{-\alpha\beta p_0+\alpha\beta+\beta-1}\,.
$$
Since
$$
II(s)\leq s^{\alpha\beta} s^{-\alpha\beta p_0-1}
s^{\alpha\beta/2}\int_0^{s^\beta} z^{-\alpha/2} dz
= c s^{-\alpha\beta p_0+\alpha\beta+\beta-1}\,,
$$
we get $II(s)\leq c I(s)$.
\qed

The following result is an immediate consequence of Lemma~\ref{ler}.
\begin{lem}
  \label{le1}
If $s\geq 1$, $x=(x_1,\tx)\in \Rd$ and $x_1\leq s/2$, then
\begin{equation}
  \label{ee1}
\PP_x\{X_\ts\in \Pb\}\leq c s^{-\alpha\beta (p_0-1)}\,.
\end{equation}
\end{lem}

Estimate (\ref{ee1}) is sharp if, for example,  $x=(s/2,\tilde{0})$, where
$\tilde{0}=(0,\ldots,0)\in \Rdl$. Indeed, let $B=B(x,cs^\beta)\subset
\Pb^{0,s}$, where  $c>0$ is small enough. By (\ref{Poisson}) we have 
\begin{eqnarray}
\PP_x\{X_\ts\in \Pb\}&\geq& 
\PP_x\{X_{\tau_B}\in \Pb^{s,\infty}\}\nonumber\\
&\geq& c(s^\beta)^\alpha\int_{s}^\infty dz_1 
\int_{\{\tz\in\Rdl\,:\; |\tz|<z_1^\beta\}} |(z_1,\tz)-x|^{-d-\alpha}
d\tz\nonumber \\ 
&\geq& cs^{\alpha\beta}\int_{s}^\infty t^{\beta(d-1)-d-\alpha} dt
=c s^{-\alpha\beta (p_0-1)}\,.
\label{eis}
\end{eqnarray}
We can, however, improve the estimate when $x$ is small relative to $s$. 
\begin{lem}
  \label{lomh}
If $s\geq 1$, $x=(x_1,\tx)\in \Rd$ and $x_1\leq s/2$, then
\begin{equation}
  \label{eomh}
\PP_x\{X_\ts\in \Pb\}\leq C (x_1\vee 1)^{\alpha \beta} s^{-\alpha\beta p_0}\,.
\end{equation}
\end{lem}
\proof
If $s\leq S$, where $0<S<\infty$ is fixed, then (\ref{eomh}) trivially
holds with $C=S^{\alpha \beta p_0}$, hence from now on we assume
$s\geq S$, where $S=2^{k_0}$ and $k_0\geq 2$ is such that
$$
\frac{4^{\alpha\beta}c}{\alpha\beta(p_0-1)}
\left(2^{k_0}\right)^{-\alpha\beta(p_0-1)}<\frac{1}{2}\,.
$$
The reason for this choice of $k_0$ will be made clear later on.
Here and in what follows $c$ is the constant of Lemma~\ref{ler}.
We will prove (\ref{eomh}) by induction.
If $s/4\leq x_1\leq s/2$ then by Lemma~\ref{le1}
\begin{eqnarray*}
\PP_x\{X_\ts\in \Pb\}&\leq& c_1 s^{\alpha\beta}s^{-\alpha\beta p_0}
\leq 4^{\alpha\beta}c_1 x_1^{\alpha\beta}s^{-\alpha\beta p_0}
=c_1 (x_1\vee 1)^{\alpha\beta}s^{-\alpha\beta p_0}\,.
\end{eqnarray*}
Assume that $n$ is a natural number and (\ref{eomh}) holds for all
$x=(x_1,\tx)\in \Pb$ such that $s/2^{n+1}\leq x_1\leq s/2$.

Let $s/2^n\geq 1$ and $x=(x_1,\tx)\in \Pb$ be such that $s/2^{n+2}\leq
x_1<s/2^{n+1}$. Note that $1\leq s/2^n\leq 4x_1$ here. We have
\begin{eqnarray*}
\PP_x\{X_\ts\in \Pb\}&\leq& 
 \PP_x\{X_{\tau(s/2^n)}\in \Pb^{s/2,\infty}\}\\
&+&\EE_x\{X_{\tau(s/2^n)}\in \Pb^{s/2^n,s/2}\,;\;
\PP_{X_{\tau(s/2^n)}}\{X_\ts\in \Pb\}\}\\
&=&I+II\,.
\end{eqnarray*}
By Lemma~\ref{ler},
\begin{eqnarray*}
I&\leq&c(s/2^n)^{\alpha\beta}\int_{s/2}^\infty t^{-\alpha\beta p_0-1}dt
=c2^{\alpha\beta p_0}(s/2^n)^{\alpha\beta}s^{-\alpha\beta p_0}\\
&\leq&c2^{\alpha\beta p_0+2\alpha}x_1^{\alpha\beta}s^{-\alpha\beta p_0}
\leq c_2(x_1\vee 1)^{\alpha\beta}s^{-\alpha\beta p_0}\,.
\end{eqnarray*}
By Lemma~\ref{ler} and induction
\begin{eqnarray*}
II&\leq&c(s/2^n)^{\alpha\beta}\int_{s/2^n}^{s/2} t^{-\alpha\beta p_0-1}
C t^{\alpha\beta}s^{-\alpha\beta p_0}dt\\
&\leq&\frac{c}{\alpha\beta(p_0-1)}(s/2^n)^{-\alpha\beta(p_0-1)}
C(s/2^n)^{\alpha\beta}s^{-\alpha\beta p_0}\\
&\leq&\frac{4^{\alpha\beta}c}{\alpha\beta(p_0-1)}(s/2^n)^{-\alpha\beta(p_0-1)}
C(x_1\vee 1)^{\alpha\beta}s^{-\alpha\beta p_0}\,.
\end{eqnarray*}
Thus
$$
R=\frac{I+II}{(x_1\vee 1)^{\alpha \beta} s^{-\alpha\beta p_0}}
\leq c_2+\frac{4^{\alpha\beta}c}{\alpha\beta(p_0-1)}(s/2^n)^{-\alpha\beta(p_0-1)}C\,.
$$
For $n$ such that $s/2^n\geq 2^{k_0}$ we have $R\le c_2+C/2$ and we
can take $C=2(c_1\vee c_2)$ in our inductive assumption to the effect that $R\leq
c_2+c_1\vee c_2\leq C$, and so (\ref{eomh}) holds for every $x=(x_1,\tx)\in \Pb$,
satisfying $s/2^{n+2}\leq x_1\leq s/2$.

By induction, (\ref{eomh}) is true with $C=2(c_1\vee c_2)$ for all
$x=(x_1,\tx)\in \Pb$, satisfying $2^{k_0-2}\leq x_1\leq s/2$.

If $x=(x_1,\tx)\in \Pb$ and  $0< x_1\leq 2^{k_0-2}$, then 
\begin{eqnarray*}
\PP_x\{X_\ts\in \Pb\}&\leq& 
\PP_x\{X_{\tau(2^{k_0}/2)}\in \Pb^{s/2,\infty}\}\\
&+&\EE_x\{X_{\tau(2^{k_0}/2)}\in \Pb^{2^{k_0}/2,s/2}\,;\;
\PP_{X_{\tau(2^{k_0}/2)}}\{X_\ts\in \Pb\}\}\\
&\leq&c(2^{k_0}/2)^{\alpha\beta}\int_{s/2}^\infty t^{-\alpha\beta p_0-1}dt\\
&+&
c(2^{k_0}/2)^{\alpha\beta}\int_{2^{k_0}/2}^{s/2} t^{-\alpha\beta p_0-1}
2(c_1\vee c_2) t^{\alpha\beta}s^{-\alpha\beta p_0}dt\\
&\leq& c_3 s^{-\alpha\beta p_0}
\leq  c_3 (x_1\vee 1)^{\alpha\beta}s^{-\alpha\beta (p_0-1)}\,.
\end{eqnarray*}
We used here our previous estimates. The proof is complete.
\qed

Note that (\ref{eomh}) is sharp if, for example, $x=(x_1,\tilde{0})$ and
$1\leq x_1\leq s/2$.  This can be verified as in (\ref{eis}) to the
effect that for such $x$ the upper bound in Lemma~\ref{lomh}
is of the same order as $\PP_x\{X_{\tau_B}\in \Pb^{s,\infty}\}$,
where $B$ is the largest ball centered at $x$ such that $B\subset
\Pb$. 

Informally, $\{X_t\}$ goes to $\Pb^{s,\infty}$ ``mostly'' by a direct jump
from $B$. 

This informal rule seems to be related to a ``thinness'' of
$\Pb$ at infinity (or its inversion \cite{BZ} at $0$).  This is false
for cones (\cite{BB}).
 
\begin{lem}\label{lB}
  If $p< p_0$, then $E_x \tb^p<\infty$ for every $x\in \Pb$.
\end{lem}
\proof
Let $1/\lambda_1$ be the first eigenvalue of the Green operator 
$G_B f(\tx)={\mathcal E}_{\tx} \int_0^{\eta_B}f(Y_t)dt$, $\tx\in \Rdl$, 
for our isotropic stable process $Y_t$  in $\Rdl$. 
Here, $B$  is the unit ball in $\Rdl$, $\eta_D$ is the first exit
time  of $Y$ from $D\subset \Rdl$, and $\mathcal P$, $\mathcal E$, 
are, respectively, the distribution  and expectation corresponding to $Y$. 
For $r>0$, by {\it scaling}, $1/(r^{-\alpha} \lambda_1)$ is the
eigenvalue for $G_{rB}$ and
$\mathcal{P}_{\tx}\{\eta_{rB}>t\}\leq c e^{-tr^{-\alpha}\lambda_1}$,
$0<t<\infty$.
Let $\mathcal C$ be the cylinder as in (\ref{edc}). For $t>0$, fixed $x=(x_1,\tx)\in \Pb$ 
and $s\geq 1\vee 2x_1$, we have by Lemma~\ref{lomh} that
\begin{eqnarray*}
  \PP_x\{\tb>t\}&=&
  \PP_x\{\tb>t\,,\; \tb=\ts\}+\PP_x\{\tb>t\,,\;\tb>\ts\}\\
&\leq& \PP_x\{\tau_{\mathcal C}>t\}
+\PP_x\{X_\ts\in \Pb\}
\leq c e^{-ts^{-\alpha\beta}\lambda_1}+c s^{-\alpha\beta p_0}\,.
\end{eqnarray*}
Let us take $s=t^{(1-\epsilon)/(\alpha\beta)}$, where
$0<\epsilon\leq 1/2$. We get
\begin{eqnarray*}
  \PP_x\{\tb>t\}
\leq c e^{-t^{\epsilon}\lambda_1}+c t^{-(1-\epsilon) p_0}\leq c t^{-(1-\epsilon)p_0}
\end{eqnarray*}
for large $t$.
Thus, letting $p=(1-2\epsilon)p_0$ we get
$$
\EE_x\tb^p=\int_0^\infty p t^{p-1}\PP_x\{\tb>t\}\,dt
\leq const.+const.\int_1^\infty t^{-\epsilon p_0-1}dt<\infty\,.
$$
This finishes the proof.
\qed

\noindent
{\it Proof of Theorem~\ref{iet}}.\/ 
The result follows from Lemma~\ref{lA} and Lemma~\ref{lB}.
\qed

\noindent
We conclude with three remarks.

Because of scaling of $\{X_t\}$, Theorem~\ref{iet} holds with the
same $p_0$ for the more general parabola-shaped regions of the form 
$$
\{x=(x_1,\tilde{x})\;:\; x_1>0\,,\; \tilde{x}\in \Rdl\,,|\tilde{x}|<a\,x_1^\beta\}, 
$$
for any  $0<a<\infty$.

If $\beta\downarrow 0$ then $\Pb$ approaches an infinite cylinder, for
which the exit time
has exponential moments (compare the proof of Lemma~\ref{lB}).

The second endpoint $\beta\uparrow 1$ suggests studying the rate of
convergence (to $1$) of the critical exponent of integrability of the exit time 
from the right circular cone as the opening of the cone tends to $0$.
A partial result in this direction is given in \cite{K3}. \cite{BB}
contains more information on our stable processes in cones 
and a hint how to approach the problem.

{\bf Acknowledgment}. The second author is grateful for the
hospitality of the Departments of Statistics and Mathematics at Purdue
University, where the paper was prepared in part.  We are grateful to the referee for pointing
out several corrections.  

\vspace{22pt}

\end{document}